\documentclass[11pt,reqno]{amsart}
\usepackage{amsfonts}
\usepackage{mathrsfs}
\usepackage{graphicx} 
\usepackage{epsfig}
\usepackage{amsmath, amsthm}
\usepackage{amssymb}

\theoremstyle{plain}
\newtheorem{thm}{Theorem}[section]
\newtheorem{prop}[thm]{Proposition}
\newtheorem{lem}[thm]{Lemma}
\newtheorem{cor}[thm]{Corollary}

\theoremstyle{definition}

\newtheorem{rem}[thm]{Remark}
\newtheorem{defn}[thm]{Definition}
\newtheorem{eg}[thm]{Example}
\newtheorem{subtitle}[thm]{}
\newtheorem{ex}{Exercise}[section]
\numberwithin{equation}{section}

\def\d{\delta}

\def\g{\gamma}

\def\K{\nabla}
\def\l{\lambda}

\def\n{\,\vert\,}
\def\o{\theta}

\def\co{{\mathcal{O}}}

\def\li{\langle}
\def\ri{\rangle}
\def\n{\ \vert\ }

\def\ms{\medskip}

\def\ni{\noindent}
\def\ti{\tilde}
\def\p{\partial}

\def\Im{{\rm Im\/}}
\def\I{{\rm I\/}}

\def\rd{{\rm d\/}}

\def\R{\mathbb{R} }
\def\C{\mathbb{C}}

\newcommand{\beg}{\begin{eg}}
\newcommand{\eeg}{\end{eg}}
\newcommand{\beq}{\begin{equation}}
\newcommand{\eeq}{\end{equation}}
\newcommand{\bthm}{\begin{thm}}
\newcommand{\ethm}{\end{thm}}
\newcommand{\bprop}{\begin{prop}}
\newcommand{\eprop}{\end{prop}}
\newcommand{\bcor}{\begin{cor}}
\newcommand{\ecor}{\end{cor}}
\newcommand{\blem}{\begin{lem}}
\newcommand{\elem}{\end{lem}}
\newcommand{\bca}{\begin{cases}}
\newcommand{\eca}{\end{cases}}
\newcommand{\brem}{\begin{rem}}
\newcommand{\erem}{\end{rem}}
\newcommand{\bpm}{\begin{pmatrix}}
\newcommand{\epm}{\end{pmatrix}}
\newcommand{\bbm}{\begin{bmatrix}}
\newcommand{\ebm}{\end{bmatrix}}
\newcommand{\bvm}{\begin{vmatrix}}
\newcommand{\evm}{\end{vmatrix}}
\newcommand{\bdefn}{\begin{defn}}
\newcommand{\edefn}{\end{defn}}
\newcommand{\bsub}{\begin{subtitle}}
\newcommand{\esub}{\end{subtitle}}
\newcommand{\bex}{\begin{ex}}
\newcommand{\eex}{\end{ex}}
\newcommand{\ben}{\begin{enumerate}}
\newcommand{\een}{\end{enumerate}}

\def\calM{\mathcal{M}}

\def\R{\mathbb{R}}
\def\calS{\mathcal{S}}

\def\R{\mathbb{R}}

\def\C{\mathbb{C}}

\def\det{{\rm det \/ }}

\def\rd{{\rm\, d\/}}

\def\fb{\mathfrak{B}}

\def\sech{{\rm sech\/}}

\def\an1{$A^{(1)}_{n-1}$}
\def\Ker{{\rm Ker\/}}
\def\rs{\R^2\setminus\{0\}}
\def\rsn{\R^n\setminus\{0\}}

\begin{document}

\title[Central Affine Curve Flow on the Plane]
{Central Affine Curve Flow on the Plane} 


\author{Chuu-Lian Terng$^\dag$}\thanks{$^\dag$Research supported
in  part by NSF Grant DMS-1109342}
\address{Department of Mathematics\\
University of California at Irvine, Irvine, CA 92697-3875.  Email: cterng@math.uci.edu}
\author{Zhiwei Wu}
\address{ Department of Mathematics\\ Ningbo University\\ Ningbo, Zhejiang, 315211, China. Email: wuzhiwei@nbu.edu.cn}


\maketitle

\centerline{\it Dedicated to Prof. Choquet-Bruhat}

\begin{abstract}

We give the following results for  Pinkall's central affine curve flow on the plane: (i) a systematic and simple way to construct the known higher commuting curve flows, conservation laws, and a bi-Hamiltonian structure, (ii)  B\"acklund transformations and a permutability formula, (iii) infinitely many families of explicit solutions. We also solve the Cauchy problem for periodic initial data. 

\end{abstract}

\section{Introduction}\label{dn}

The group $SL(2,\R)$ acts transitively on $\R^2\setminus\{0\}$ by $A\cdot y= Ay$.  It is noted in \cite{UP95} that given a smooth curve $\g$ on $\rs$, if $\det(\g, \g_s)$ never vanishes, then there is a unique parameter $x$ such that
$$\det(\g, \g_x)=1.$$ (In fact, $\frac{\rd s}{\rd x}= \det(\g, \g_s)^{-1}$).  Taking $x$-derivative of $\det(\g, \g_x)=1$ gives $\det(\g, \g_{xx})=0$. Hence there is a unique smooth function $q$ such that
$$\g_{xx}= q\g.$$ This parameter $x$ is called the {\it central affine arc-length parameter\/} and $q$ is called the {\it central affine curvature\/} of $\g$. Note that $q=\det(\g_{xx}, \g_x)$.

Let
\beq\label{do}
\calM_2(I)=\{\g:I\to \rs \n \g \,\, {\rm smooth\, curve,\,} \det(\g, \g_x)=1\},
\eeq
where $\I$ is $\R$ or $S^1$. 

Note that $X$ lies in the tangent space  $T(\calM_2(I))_\g$ of $\calM_2(I)$ at $\g$ if and only if $\det(X, \g_x)+ \det(\g, X_x)=0$.  So $X=y_1\g+ y_2\g_x$ lies in $T(\calM_2(I))_\g$ if and only if
$$y_1= -\frac{(y_2)_x}{2}.$$
This identifies $T(\calM_2(I))_\g$ as $C^\infty(I, \R)$.
Henceforth we will use the following notation: given $\xi\in C^\infty(I, \R)$, let $\ti \xi$ denote the tangent vector field on $\calM_2(I)$ defined by
\beq\label{dt}
\ti\xi(\g)= -\frac{\xi_x}{2}\g + \xi \g_x.
\eeq

We call the flow on $\calM_2(I)$ defined by a tangent vector field $\ti\xi(\g)$ on $\calM_2(I)$,
$$\g_t= \ti\xi(\g) =-\frac{\xi_x}{2} \g + \xi \g_x,$$
 a {\it central  affine curve flow on $\rs$\/} if $\xi$ is a differential polynomial of the central affine curvature $q$ for $\g$ (i.e., $\xi$ is a polynomial of $q$ and the $x$-derivatives of $q$).  
 Note that a central affine curve flow is invariant under the action of $SL(2,\R)$ and translations, i.e.,
if $\g$ is a solution of \eqref{cf}, $c\in SL(2,\R)$, and $x_0, t_0$ are real constants, then
$c\g$ and $\g_1(x, t)=  \g(x+x_0, t+ t_0)$ are again solutions of \eqref{cf}.

Pinkall considered the following in \cite{UP95} third oder central affine curve flow on $\calM_2(I)$:
\begin{equation}\label{cf}
\gamma_t=\frac{1}{4}q_x \gamma-\frac{1}{2}q\gamma_x,
\end{equation}
i.e., the flow defined by the vector field $\ti \xi$ with $\xi=
-\frac{q}{2}$. He  proved in \cite{UP95} that if $\gamma$ is a solution of the curve
flow \eqref{cf}, then the {\it central affine curvature} $q(\cdot, t)$ of $\g(\cdot, t)$
satisfies the KdV equation:
\begin{equation}\label{kdv}
q_t=\frac{1}{4}(q_{xxx}-6qq_x).
\end{equation}
He also proved in \cite{UP95} that 
\beq\label{dp}
w_\g(\ti\xi, \ti\eta)= \oint \det(\ti\xi, \ti\eta)=\frac{1}{2} \oint (\eta_x \xi-\xi_x \eta) \rd x =-\oint \xi_x \eta\rd x,
\eeq
is a symplectic form on the orbit space $\calM_2(S^1)/S^1$ and \eqref{cf} is the Hamiltonian flow for 
$$H(q)= \oint \frac{1}{2} q\rd x.$$ Here $S^1$ acts on $\calM_2(S^1)$ by 
$$(e^{i\o}\cdot \g)(x)= \g(x+\o).$$

Chou and Qu  wrote down the traveling wave solutions for the central affine curve flow \eqref{cf} in \cite{CQ01}.  Calini, Ivey and Mari-Beffa studied periodic (in $x$) solutions of the curve flow \eqref{cf} whose central affine curvatures are finite-gap solutions of the KdV equation in \cite{CIM09}. Higher order central affine curve flows and conservation laws for \eqref{cf} were given in \cite{CQ01}, \cite{CIM09},  \cite{FK10}. A bi-Hamiltonian structure for \eqref{cf} was discussed in \cite{FK13}.

In this paper, 
\ben
\item We explain how to use soliton theory of the KdV hierarchy to give a systematic and simple way to obtain higher commuting Hamiltonian flows, conservation laws, and a bi-Hamiltonian structure for the central affine curve flow \eqref{cf}.  

\item We construct B\"acklund transformations and a Permutability formula for the central affine curve flow \eqref{cf}. Then we use these to construct recursively, infinitely many families of explicit solutions of \eqref{cf} whose central affine curvatures are pure 
n-soliton solutions or rational solutions of the KdV equation.

\item We solve the Cauchy problem for \eqref{cf} with initial data $\g_0$ with rapidly decaying central affine curvature. We also solve the periodic Cauchy problem for \eqref{cf}.
\een

A natural generalization of the KdV hierarchy to $(n-1)$ component functions is
the {\it \an1-KdV hierarchy \/}, a KdV type hierarchy
constructed by Drinfel'd and Sokolov from the affine Kac-Moody
algebra \an1 in \cite{DirSok84}. This inspired the study of
$n$-dimensional central affine curve flows: It follows from the change
of variable formula that if $\g$ is a smooth curve in $\rsn$ such
that
$$\det(\g, \g_s, \ldots, \g_s^{(n-1)})>0,$$
then there is a unique orientation preserving parameter $x$ such that
$$\det(\g, \g_x, \ldots, \g_x^{(n-1)})\equiv 1.$$
Taking the $x$-derivative of the above equation gives
$$\det(\g, \g_x, \ldots, \g_x^{(n-2)}, \g_x^{(n)})\equiv 0.$$
So there exist unique smooth functions $u_1, \ldots, u_{n-1}$ such that
$$\g_x^{(n)}= u_1 \g + u_2 \g_x + \cdots + u_{n-1} \g_x^{(n-1)}.$$
This parameter $x$ is called the {\it central affine arc-length
parameter\/} and $u_i$ is called the {\it $i$-th central affine
curvature\/} of $\g$ for $1\leq i\leq n-1$.  It follows from the
existence and uniqueness for ordinary differential equations that these $u_i$ form a complete
set of local invariants for curves in $\rsn$ under the group
$SL(n,\R)$. In a forthcoming paper \cite{TW2}, we consider the
following central affine curve flow \beq\label{ncf2}
\gamma_t=-\frac{2}{n}u_{n-1}\gamma+\gamma_{xx}, \eeq 
on
$$\calM_n(I)=\{\g:I\to \R^n\setminus\{0\}\n \det(\g, \g_x, \ldots,
\g_x^{(n-1)})=1\},$$ where $u_{n-1}(\cdot, t)$  is the $(n-1)$-th
central affine curvature for $\gamma(\cdot, t)$.  When $n=3$, this
central affine curvature flow was studied in \cite{HS02} and
\cite{CIM13}. In the forthcoming paper
\cite{TW2}, we (i) prove that if $\g(x,t)$ is a solution
of \eqref{ncf2} on $\calM_n(I)$ then its central affine curvatures $u_1(\cdot, t),
\ldots, u_{n-1}(\cdot, t)$ satisfy the second flow in the \an1-KdV
hierarchy, (ii) construct for the curve flow \eqref{ncf2} higher order commuting central affine
curve flows, conservation laws, a bi-Hamiltonian structure, B\"acklund transformations, and (iii) obtain  recursively, infinitely many families of explicit solutions of \eqref{ncf2}. 

This paper is organized as follow: In section \ref{bq}, we construct a sequence of commuting higher order central affine curve flows on $\rs$ and solve the Cauchy problem. In section \ref{dv}, we construct B\"acklund transformations and a
Permutability formula for the central affine curve flow \eqref{cf}
and then apply these to the stationary solution of \eqref{cf} to
obtain infinitely many families of explicit solutions of \eqref{cf} whose
affine curvatures are pure $n$-soliton solutions or rational
solutions of the KdV equation.  We discuss the
Hamiltonian aspsect of \eqref{cf} in the final section.

\bigskip
\section{Higher order commuting curve flows and Cauchy problem}\label{bq}

The outline of this section is as follows:
\ben
\item We review the construction of the KdV hierarchy and its Lax pair (for detail see \cite{AKNS74}, \cite{Dic03}, \cite{TU11}).
\item We use the Lax pairs of the KdV hierarchy to write down the known sequence of commuting higher order central affine curve flows on $\calM_2(I)$ constructed in \cite{CQ01}, \cite{CIM09}, \cite{FK10},  whose central affine curvatures satisfy the higher flows in the KdV hierarchy. 
\item We use the solution of the Cauchy problem for the KdV equation to solve the Cauchy problem for \eqref{cf} with initial data $\g_0\in\calM_2(\R)$ having rapidly decaying central affine curvatures. We also solve the Cauchy problem for \eqref{cf} with periodic initial data $\g_0\in \calM_2(S^1)$.
\een

\ms
Let $\fb:sl(2,\R)\to \R e_{12}$ denote the linear map defined by
\beq\label{ao}\fb \left(\bpm  a&b\\ c&-a\epm\right)=\bpm 0&c\\ 0&0\epm.
\eeq

Given a smooth function $q:\R\to \R$, it is known (cf. \cite{TU11}) that there exists a unique
$$Q(q,\l)= e_{12}\l  + \sum_{i\geq 0} Q_{-i}(q)\l^{-i}$$
satisfying
\beq\label{ar}
\bca \left[\rd_x + \bpm 0& \l + q\\ 1&0\epm, Q(q,\l)\right]=0,\\ Q(q,\l)^2= \l\I_2.\eca
\eeq
Compare coefficients of $\l^{-j}$ of the above equation to get the following recursive formulas:
\beq\label{du}
\bca (Q_{-j}(q))_x+[e_{21} + qe_{12}, Q_{-j}(q)] =[Q_{-(j+1)}(q), e_{12}],\\
e_{12}Q_{-(j+1)}(q)+ Q_{-(j+1)}(q)e_{12} +\sum_{i=0}^j Q_{-i}(q) Q_{-(j-i)}(q)=0,\eca
\eeq
for all $j\geq 0$. Write
$$Q_{-j}(q)=\bpm A_j(q) & B_j(q)\\ C_j(q) & -A_j(q)\epm.$$
Compare entries of \eqref{du} to get
\begin{align}
&C_{j+1}(q)= -((A_j(q))_x + qC_j(q)-B_j(q)), \label{dq1}\\
& A_{j+1}(q)= \frac{1}{2}((B_j(q))_x - q A_j(q),\label{dq2}\\
& A_j(q)= -\frac{1}{2} (C_j(q))_x. \label{dq3}
\end{align}
It follows that these $Q_i(q)$'s can be obtained recursively and they are differential polynomials of $q$ in $x$-variable.
For example, 
\begin{align*}
&Q_0(q)= \bpm 0& \frac{q}{2}\\ 1&0\epm, \quad Q_{-1}(q)=\frac{1}{4} \bpm q_x& q_{xx}-2q^2\\ -2q& -q_x\epm, \\
&Q_{-2}(q)=\frac{1}{8} \bpm \frac{1}{2}q_{xxx}-3qq_x & \frac{1}{4}(q^{(4)}-6qq_{xx}-7q_x^2+2q^3)\\ 3q^2 -q_{xx} &-\frac{1}{2}q_{xxx}+3qq_x\epm, \\
&Q_{-3}(q)= \bpm \ast &  \ast \\
-\frac{1}{32}(q^{(4)}+10q^3-5q_x^2-10qq_{xx}) & \ast \epm.
\end{align*}
{\it The $(2j+1)$-th flow in the KdV hierarchy\/} is
\beq\label{aq}
q_{t_{2j+1}}= (B_j(q)-C_{j+1}(q))_x - 2q A_j(q),
\eeq
In particular, the first, third, and fifth flows are:
\begin{align}
&q_{t_1}= q_x, \label{bi1}\\
&q_{t_3}= \frac{1}{4}(q_{xxx} - 6 qq_x), \label{bi3}\\
&q_{t_5}=\frac{1}{16}(\p_x^5 q-10\,q\p_x^3 q-20(\p_x q)\p_{xx}q+30\,q^2\p_xq). \label{bi5}
\end{align}
Note that the third flow \eqref{bi3} is the KdV equation.

The following two Theorems are well-known and the proofs can be found in many places (cf. \cite{ZF71}, \cite{AKNS74}, \cite{TU11}).

\bthm  The flows in the KdV hierarchy commute.
\ethm

\bthm \label{at} {\rm [Lax pair for the KdV hierarchy]} \

 \ni The following statements are equivalent for  $q\in C^\infty(\R^2,\R)$:
 \ben
 \item[(i)] $q$  is a solution of the $(2j+1)$-th flow \eqref{aq}.
\item[(ii)] The following family of connections on the $(x, t_{2j+1})$-plane defined by $q$ is flat for all parameter $\l\in \C$,
 \beq\label{am}
\left[\p_x+ \bpm 0& q+\l\\ 1&0\epm,\, \p_{t_{2j+1}}+ (Q(q,\l)\l^j)_+ -\fb(Q_{-(j+1)}(q))\right]=0,
\eeq
where $\fb$ is linear map defined by \eqref{ao}. We call \eqref{am} the {\it Lax pair\/} of the solution $q$ of the $(2j+1)$-th flow in the KdV hierarchy. 
\item[(iii)] Equation \eqref{am} holds for $\l=0$, i.e.,
\beq\label{an} \left[\p_x+ \bpm 0& q\\ 1&0\epm, \, \p_{t_{2j+1}} +
Q_{-j}(q)- \fb(Q_{-(j+1)}(q))\right]=0. \eeq \een \ethm

\ms

For example, \eqref{am} for the third flow (KdV) is 
\beq\label{ap}
\left[\p_x+ \bpm 0& \l+q\\ 1&0\epm, \, \p_{t_{3}} + \bpm
\frac{1}{4}q_x &\l^2+\frac{1}{2} q\l +\frac{1}{4}(q_{xx}- 2q^2) \\
\l-\frac{q}{2} & -\frac{1}{4}q_x\epm\right]=0. \eeq

Let $q$ be a solution of \eqref{aq}, and $c(\l)\in SL(2,\C)$ holomorphic for $\l\in \C$ satisfying $\overline{c(\bar\l)}= c(\l)$. Then there exists a unique $E(x,t,\l)\in SL(2,\C)$ satisfying $\overline{E(x,t,\bar\l)}= E(x,t,\l)$ and
$$\bca E_x= E\bpm 0& \l+q\\ 1&0\epm,\\
E_t= E((Q(q,\l)\l^j)_+ -\fb(Q_{-(j+1)}(q))), \\ E(0,0, \l)= c(\l).
\eca$$ Moreover, the solution $E(x,t,\l)$ is holomorphic for $\l\in
\C$.  We call $E$ an
{\it extended frame\/} of the solution $q$ of the $(2j+1)$-th flow \eqref{aq}.

\ms

Next we discuss Pinkall's result that the central affine curvature $q(\cdot, t)$ of a solution of \eqref{cf} is a solution of the KdV equation.
Let $\g$ be a solution of \eqref{cf} on $\calM_2(I)$, and $g=(\g, \g_x)$. Then we have $\g_{xx}=q \g$ and $g_x= g\bpm 0& q\\ 1&0\epm$.  Use $\g_t= \frac{1}{4} q_x\g -\frac{1}{2} q\g_x$ to get
$$(\g_x)_t= (\g_t)_x=\left( \frac{1}{4} q_x\g -\frac{1}{2} q\g_x\right)_x= (\frac{1}{4} q_{xx}- \frac{1}{2} q^2) \g -\frac{1}{4} q_x \g_x.$$
So $g=(\g, \g_x)\in SL(2,\R)$ satisfies
\beq\label{ac}
\bca g_x= g\bpm 0& q\\ 1&0\epm,\\ g_t= g \bpm \frac{1}{4} q_x & \frac{1}{4}(q_{xx}- 2q^2)\\ -\frac{q}{2}& -\frac{1}{4}q_x\epm. \eca
\eeq
This implies that
$$\left[\p_x+\bpm 0& q\\ 1&0\epm, \, \p_t+ \bpm \frac{1}{4} q_x & \frac{1}{4}(q_{xx}- 2q^2)\\ -\frac{q}{2}& -\frac{1}{4}q_x\epm\right]=0.$$ By Theorem \ref{at},  $q$ is a solution of the KdV equation. This gives Pinkall's result:

\bprop \label{ah} ( \cite{UP95}) If  $\gamma$ is a solution of
\eqref{cf} on $\calM_2(I)$, then its central affine curvature $q(\cdot,
t)$ of $\gamma(\cdot, t)$ is a solution of the KdV equation
\eqref{kdv}. \eprop

The converse is also true when $\I=\R$:

\bprop\label{cfkdv}
If $q:\R^2\to \R$ is a solution of the KdV equation and $c_0\in SL(2,\R)$, then
\ben
\item[(i)] there is a unique $g:\R^2\to SL(2,\R)$ satisfies \eqref{ac} with $g(0,0)= c_0$,
\item[(ii)] $\g(x, t)=g(x,t)p_0$ is a solution of the
curve flow \eqref{cf}  with central affine curvature $q$, where $p_0=(1,0)^t$.
\een
\eprop

\begin{proof} Statement (i) follows from properties of the Lax pair.
Let $\g$ and $v$ denote the first and the second column of $g$ respectively. Then the first equation of \eqref{ac} implies that $\g_x= v$ and $v_x= q\g$. So $\g_{xx}= q\g$. Since $g=(\g, \g_x)\in SL(2,\R)$, $x$ is the central affine arc-length for the curve $\g(\cdot, t)$ and $q$ is its central affine curvature. The second equation of \eqref{ac} implies that
$\gamma_t=\frac{1}{4}q_x\gamma-\frac{1}{2}q\gamma_x$.
\end{proof}

If $g_0$ is the solution of \eqref{ac} with $g_0(0,0)=\I_2$, then given $c\in SL(2,\R)$ the solution $g$ of \eqref{ac} with $g(0,0)= c$ is $cg_0$.  So it follows from Proposition \ref{cfkdv} that we have

\bcor\label{bj} Let $\Psi:\calM_2(I)\to C^\infty(I, \R)$ be the map defined by $\Psi(\g)=$ the central affine curvature of $\g$. Then
 $\Psi(\g_1)=\Psi(\g_2)$ if and only if there is a constant $c\in SL(2,\R)$ such that $\g_2= c\g_1$. 
 \ecor
 
 We define the holonomy map next:
 
 \bdefn Given $q\in C^\infty(S^1, \R)$, The {\it holonomy map\/} $\Pi$ is the map from $ C^\infty(S^1,\R)$ to $ SL(2,\R)$ defined by $\Pi(q)= $ the holonomy of the connection $\frac{\rd}{\rd x}+ \bpm 0& q\\ 1&0\epm$, i.e., $\Pi(q)= g(2\pi)$, where $g$ is the solution of 
 $$\bca g^{-1}g_x= \bpm 0& q\\ 1 & 0\epm,\\ g(0)=\I_2.\eca$$
 \edefn
 
 \bcor \label{eq} Let $\Psi:\calM_2(S^1)\to C^\infty(S^1,\R)$ be as in Corollary \ref{bj} and $\Pi$ the holonomy map defined above. Then $\Psi$ induces a bijection from the orbit space $\calM_2(S^1)/SL(2,\R)$ onto 
 $$C^\infty_I(S^1,\R)=\{q\in C^\infty(S^1, \R)\n \Pi(q)= \I_2\}.$$ 
 \ecor

\ms Next we write down Commuting higher order central affine curve flows for \eqref{cf}.  Recall that $Y(\g)= y_1\g + y_2\g_x$ is a tangent vector field
on $\calM_2(I)$ if and only if $y_1=-\frac{1}{2} (y_2)_x$, and we use
$\ti y_1$ to denote this vector field.   So it follows from
\eqref{dq3} that  $A_j(q)\g + C_j(q)\g_x$ is tangent to $\calM_2(I)$ at
$\g$ and 
\beq\label{cfj} 
\g_{t_{2j+1}}= A_j(q) \g+ C_j(q)\g_x= -\frac{1}{2}
(C_j(q))_x\g + C_j(q)\g_x \eeq is a central affine curve flow on
$\calM_2(I)$ of order $2j+1$, where $Q_{-j}(q)= \bpm A_j(q) & B_j(q)\\ C_j(q)
& -A_j(q)\epm$ is the coefficient of $\l^{-j}$ of the solution
$Q(q,\l)$ of \eqref{ar}. We call this the {\it $(2j+1)$-th central
affine curve flow\/} on $\calM_2(I)$.   For example, the first, third,
and fifth (i.e., $j=0,1,2$) central affine curve flow on $\calM_2(I)$
is
\begin{align*}
&\g_{t_1}= \g_x, \\
&\g_{t_3}= \frac{1}{4} q_x \g- \frac{1}{2}q\g_x,\\
&\g_{t_5}= \frac{1}{16} (q_{xxx}- 6qq_x)\g + \frac{1}{8} (3q^2 - q_{xx}) \g_x\\
\end{align*}
Note that the third central affine curve flow is the curve flow \eqref{cf}.

\ms

The same proof of Proposition \ref{ah} implies that if $\g$ is a solution of the central affine curve flow \eqref{cfj} then its affine curvature $q(\cdot, t)$ is a solution of $(2j+1)$-th flow \eqref{aq}. Analogous result as Proposition \ref{cfkdv} for these higher order central affine curve flow can be proved in a similar manner. Since the flows in the KdV hierarchy commute, these central affine curve flows commute. Hence we get the following result proved in \cite{CQ01}, \cite{CIM09}, and \cite{FK10}:

\bprop\label{aha}
Let $\Psi$ denote the map defined in Corollary \ref{bj}. Then $\Psi$ maps the central affine curve flow \eqref{cfj} to the $(2j+1)$-th flow \eqref{aq} in the KdV hierarchy for all $j\geq 0$. Moreover, 
\ben
\item if $q$ is a solution of the $(2j+1)$-th flow of the KdV equation and $g:\R^2\to SL(2,\R)$ a solution of $g^{-1}g_x= \bpm 0& q\\ 1 & 0\epm$, then $\g= g\bpm 1\\ 0\epm$ is a solution of \eqref{cfj},
\item these flows \eqref{cfj} for all $j\geq 0$ commute. 
\een
\eprop

\ms

We will discuss  Cauchy problems for the central affine curve flow \eqref{cf} in the rest of this section.  First recall that the Cauchy problem for the KdV equation is solved for two classes of initial data $q_0$ (\cite{GGKM74}, \cite{L75}): 
\ben
\item $q_0\in \calS(\R, \R)$, i.e., $q_0$ is smooth and rapidly decaying,
\item $q_0\in C^\infty(S^1, \R)$.
\een

Use Proposition \ref{cfkdv} and the solutions to the Cauchy problem for the KdV equation with initial data $q_0\in \calS(\R,\R)$ to get

\bthm {\rm [Cauchy problem on the line]} \

\ni Suppose $\gamma_0\in \calM_2(\R)$ 
has rapidly decaying central affine curvature $q_0$, and
$q(x,t)$ the solution of the KdV equation with initial data $q(x,0)= q_0(x)$. Let $g:\R^2\to SL(2,\R)$ denote the solution of
\eqref{ac}
with  $g(0,0)=(\g_0(0), (\g_0)_x(0))$.  Then $\gamma(x,t)=g(x,t)(1,0)^t$ is the solution of \eqref{cf}  with $\gamma(x, 0)=\gamma_0(x)$. Moreover, the central affine curvatures of $\g(\cdot, t)$ are also rapidly decaying.
\ethm

If $q(x,t)$ is a solution of the KdV equation such that $q$ is periodic in
$x$ with period $2\pi$ and $g$ is a solution of \eqref{ac}, then by
Proposition \ref{cfkdv}, $g(1,0)^t$ is a solution of \eqref{cf}.
Although $q$ is periodic in $x$, the solution $g$ of the linear system
\eqref{ac} may not be periodic in $x$.  We prove below that if $g(\cdot,
0)$ is periodic then $g(\cdot, t)$ is periodic.  

\bthm {\rm [Cauchy Problem with periodic initial data]}  \label{c}\

\ni Suppose $q_0$ is the central affine curvature of $\gamma_0\in \calM_2(S^1)$ and $q(x,t)$ is the solution of the KdV equation periodic in $x$
such that $q(x,0)= q_0(x)$. Let $g:\R^2\to SL(2,\R)$ be the solution of
\eqref{ac} with initial data $g(0,0)= (\g_0(0), (\g_0)_x(0))$. Then
$\g(x,t)= g(x,t)(1,0)^t$ is a solution of \eqref{cf} with initial
data $\g(x,0)= \g_0(x)$ and $\g(x,t)$ is periodic in $x$ with period
$2\pi$. \ethm

\begin{proof} Both $g(x,0)$ and $(\g_0, (\g_0)_x)$ satisfy $g_x=g \bpm 0& q\\1 &0\epm$ with the same initial data. So by the uniqueness for ordinary differential equations we have $g(\cdot,0)= (\g_0, (\g_0)_x)$.
It follows from Proposition \ref{cfkdv} that $\gamma(x,t)=g(x,
t)(1,0)^t$ is a solution of the curve flow \eqref{cf} with initial
data $\gamma(x,0)=\gamma_0(x)$. It remains to prove that $\gamma$ is
periodic in $x$.

Since $\gamma_0$ is periodic with period $2\pi$, $g(x,0)=(\g_0(x), (\g_0)_x(x))$ is periodic in $x$. Hence $g(2\pi,0)-g(0,0)=0$. We claim that
$$y(t)= g(2\pi, t)- g(0, t)$$
is identically zero. To see this, first recall that
$$Q_{-1}(q)= \frac{1}{4}\bpm q_x & q_{xx} - 2q^2\\ -2q & -q_x\epm.$$ Since $q$ is periodic in $x$, so is $Q_{-1}(q)$.  Compute directly to get
\begin{align*}
y_t&=g_t(2\pi, t)-g_t(0, t) \\
&=(gQ_{-1}(q))(2\pi, t)-(gQ_{-1}(q)))(0, t)=(g(2\pi, t)-g(0, t))Q_{-1}(q)((0, t)\\
&=yQ_{-1}(q)(0, t).
\end{align*}
Since $y\equiv 0$ is a solution of $y_t= yQ_{-1}(q)(0,t)$ with $y(0)=0$.  But $y(0)=g(2\pi, 0)-g(0, 0)=0$. So it follows from the uniqueness for ordinary differential equations that $y$ is identically zero.
\end{proof}

\bigskip

\section{B\"acklund transformations}\label{dv}

We construct B\"acklund transformations and a Permutability formula for the central affine curve flow \eqref{cf}. Then we apply these transformations to the stationary solutions of \eqref{cf} to construct recursively, infinitely many families of explicit solutions of \eqref{cf}.

B\"acklund transformations for the KdV equation were given in several places (cf. \cite{A81}, \cite{SS93}, \cite{TU00}).  The one we will use to construct B\"acklund transformations for the central affine curve flow \eqref{cf} is given in
\cite{TU00}.  These transformations were constructed using the Lax pair for the KdV equation constructed in \cite{AKNS74}:
\beq\label{as}
\left[ \p_x + \bpm z& q\\ 1 &-z\epm, \, \p_{t_3} + \bpm z^3 -\frac{qz}{2} +\frac{ q_x }{4}& qz^2 -\frac{q_x}{2} +\frac{q_{xxx}- 2q^2}{4}\\ z^2- \frac{q}{2} & -z^3 + \frac{qz}{2} - \frac{q_x}{4} \epm\right] =0.
\eeq
This Lax pair is gauge equivalent to the Lax pair \eqref{ap} with $\l=z^2$ by 
$$\phi(z)= \bpm 1& z\\ 0&1\epm,$$ i.e.,
$$\phi(z) \left[\p_x+ \bpm 0& z^2+q\\ 1&0\epm, \, \p_{t_{3}} + \bpm
\frac{q_x}{4} &z^4+\frac{qz^2}{2} +\frac{q_{xx}- 2q^2}{4} \\
z^2-\frac{q}{2} & -\frac{q_x}{4}\epm\right]\phi(z)^{-1}$$ 
is equal to the left hand side of \eqref{as}. 
Hence $E(x,t,\l)$ is an extended frame for the Lax pair
\eqref{ap} if and only if
$$F(x,t,z)=\phi(z)E(x,t,z^2)\phi(z)^{-1}$$ is an extended frame for the Lax pair \eqref{as} in $z$-parameter.

\ms
Next we use this conjugation to state results for B\"acklund transformations in \cite{TU00}.
Given real constants $k, \xi$, we call
\beq\label{ck}
r_{\xi,k}(\l)=\bpm \xi& \xi^2-k^2+\l \\ 1& \xi\epm
\eeq
a {\it simple factor}.  A direct computation implies that
$$r_{\xi, k}^{-1}(\l)= \frac{r_{-\xi, k}(\l)}{\l-k^2}.$$
Note that $\det(r_{\xi, k}(\l))= k^2-\l$, hence $r_{\xi,k}(\l)$ has a zero at $\l=k^2$. 

\bthm\label{ab} (\cite{TU00})
 Let $q$ be a solution of the KdV equation, $E(x,t,\l)$ an extended frame of the Lax pair \eqref{ap}, and $k,\xi$ real constants.
 Set
 \begin{align*}
 &\bpm y_1(x,t)\\ y_2(x,t)\epm = E(x,t,k^2)^{-1}\bpm -\xi\\
1\epm,\\& \tilde \xi(x,t) = - \frac{y_1(x,t)}{y_2(x,t)}.
\end{align*}
If $y_2$
does not vanish in an open subset $\co\subset R^2$, then
$$r_{\xi,k}\ast q:= -q + 2(
\tilde\xi^2(x,t)-k^2)$$ is a
solution of the KdV equation defined on $\co$. Moreover,
\beq\label{ba}
 \tilde E(x,t,\l)=\frac{r_{\xi,k}(\l) E(x,t,\l)r_{-\tilde
\xi(x,t),k}(\l)}{ \l-k^2}
\eeq
 is holomorphic for $z\in \C$ and is an extended frame of
$r_{\xi,k}\ast q$.
\ethm

\brem For $k\not=0$, to prove $\ti E(x,t,z)$ defined by \eqref{ba} is holomorphic for all $\l\in \C$, it suffices to show that the residues of $\ti E$ at $\l=k^2$ are zero.
For $k=0$, we check that the constant coefficient of $r_{\xi, 0}(\l) E(x,t,\l) r_{-\ti \xi, 0}(\l)$ as a power series in $\l$ is zero.  Hence
$$\ti E(x,t,\l)= \l^{-1} r_{\xi, 0}(\l) E(x,t,\l) r_{-\ti \xi, 0}(\l)$$
is holomorphic for $z\in \C$.
\erem

\ms

Let $E$ be an extended frame of a solution $q$ of the KdV equation. Assume that
$\ti E$ is of the form \eqref{ba} and require that $\ti E^{-1}\rd
\ti E$ equals to the Lax pair of some solution $\ti q$.  This gives a system of compatible non-linear ODEs for $\ti \xi$:

\bthm\label{bm} (\cite{TU00}) Let $k\in \R$ be a constant, and
$q:\R^2\to \R$ a smooth function.  Then the following first order
system for $A$ is solvable if and only if $q$ is a solution of the
KdV equation:
$$(BT)_{q, k}\quad  \bca A_x= q-A^2+k^2,\\ A_t= \frac{q_{xx}-2q^2}{4} -\frac{q_x}{2}A + \frac{q(A^2+k^2)}{2} -k^2(A^2-k^2).
\eca$$
Moreover, if $q$ is a solution of the KdV equation and $A$ the solution of $(BT)_{q,k}$ with $A(0,0)=\xi$, then
 $A=\ti\xi$ and
 $$r_{\xi, k}\ast q:= -q + 2(\ti\xi^2-k^2)$$ is a solution of the KdV equation, where $\ti\xi$ is defined in Theorem \ref{ab}.
\ethm

\brem To construct B\"acklund transformations for a given solution $q$ of the KdV equation, we can either solve $E$ from the following linear system 
$$\bca E_x = E \bpm 0& k^2+q\\ 1&0\epm,\\
E_t=E\bpm \frac{1}{4}q_x & k^4 +\frac{1}{2} k^2 q+\frac{1}{4} (q_{xx}- 2q^2)\\ k^2- \frac{1}{2}q & -\frac{1}{4}q_x\epm,\eca$$ for constant parameter $k^2$
 or solve the non-linear system $(BT)_{q,k}$.
\erem

\ms

Permutability Theorem for the KdV hierarchy follows from a
relation among {\it simple factors\/} $p_{\xi, k}$'s:

\bprop\label{bda} (\cite{TU00}) Let $\xi_1\not= \xi_2, k_1\not= k_2$ be real constants, $\eta_1= -\xi_2
+\frac{k_1^2-k_2^2}{\xi_1-\xi_2}$, $\eta_2= -\xi_1 +
\frac{k_1^2-k_2^2}{\xi_1-\xi_2}$, and  $r_{\xi,k}$ be defined as in
\eqref{ck}.   Then $r_{\eta_2, k_2}r_{\xi_1, k_1}= r_{\eta_1,
k_1}r_{\xi_2, k_2}$. \eprop

\bthm \label{bd} {\rm [Permutability]} (\cite{TU00})\

\ni Let $\xi_1, \xi_2, k_1, k_2, \eta_1, \eta_2$ be real constants as in Proposition \ref{bda}.  Let  $E$ be an extended frame of the solution $q$ of the KdV equation, $(y_{1i}, y_{2i})^t= E(x,t, k_i^2)^{-1}(-\xi_i, 1)^t$, $\ti \xi_i=-(y_{1i}/y_{2i})$,  and
$$q_i= r_{\xi_i,k_i}\ast q = -q+ 2(\ti\xi_i^2-k_i^2), $$
for $i=1,2$.  Set
\begin{align}
&\ti\xi_{12}= -\ti\xi_1 + \frac{k_1^2-k_2^2}{\ti\xi_1-\ti\xi_2}, \label{bn1}\\
& E_{12}(x,t,\l)=\frac{r_{\eta_2, k_2}(\l) r_{\xi_1, k_1}(\l) E(x,t,\l) r_{-\ti \xi_1, k_1}(\l)r_{-\ti \xi_{12}, k_2}(\l)}{(\l-k_1^2) (\l-k_2^2)}, \label{bn2}\\
&q_{12}= -q_1+ 2(\ti\xi_{12}^2- k_2^2).\label{bn3}
\end{align}
 Then
 \ben
 \item $q_{12}= r_{\eta_2, k_2}\ast(r_{\xi_1,k_1}\ast  q)= r_{\eta_1, k_1}\ast(r_{\xi_2, k_2}\ast q)$ is a solution of the KdV equation,
 \item $E_{12}$ is an extended frame for $q_{12}$,
 \item $\ti \xi_{12}$ is the solution of ${\rm{(BT)}}_{q_1, k_2}$ with initial data $\eta_2$ and is also the solution of ${\rm {(BT)}}_{q_2, k_1}$ with initial data $\eta_1$.
  \een
\ethm

As a consequence of Proposition \ref{ah}, Proposition \ref{cfkdv}, and Theorem \ref{ab} we obtain:

\bthm\label{ax}  {\rm [BT for central affine curve flow \eqref{cf} with $k\not=0$]} \

\ni Let $\gamma$ be a solution for \eqref{cf} with central affine curvature $q$, and
$E(x, t, \l)$ the extended frame for $q$ for the Lax pair \eqref{ap} with $E(0,0,0)= (\g, \g_x)$ at $(0,0)$. Given $k, \xi\in \R$ with $k \neq 0$, let
$$\bpm y_1(x, t) \\ y_2(x, t)\epm:= E(x,
t, k^2)^{-1}\bpm -\xi \\ 1\epm.$$
Suppose $y_2(x, t)\neq 0$. Set $\ti \xi= -\frac{y_1}{y_2}$. 
Then   $$\ti \g=\frac{1}{k} (\ti\xi \gamma-\gamma_x)$$ is a solution of \eqref{cf} with central affine curvature $\ti q=r_{\xi,k}\ast q= -q+ 2(\ti \xi^2-k^2)$.  Moreover,  $$\ti E(x,t,\l)=\frac{ r_{\xi, k} (\l) E(x,t,\l) r_{-\ti\xi(x,t), k}(\l)}{\l-k^2}$$ is an extended frame for $\ti q$,
where $r_{k, \xi}(\l)=  \bpm \xi& \xi^2-k^2+\l\\ 1 & \xi\epm$.
\ethm

\begin{proof}
Since both $E(x,t,0)$ and $(\g, \g_x)$ are solutions of \eqref{ac} with the same initial data at $(0,0)$, we have $E(\cdot, \cdot,0)= (\g, \g_x)$.   By Theorem \ref{ab},
$$\ti E(x,t,\l)=\frac{ r_{\xi, k} (\l) E(x,t,\l) r_{-\ti\xi(x,t), k}(\l)}{\l-k^2}$$
 is an extended frame for the new solution $\ti q$.
By Proposition \ref{cfkdv}, the first column of $\ti E(x,t,0)$ is a solution of the curve flow \eqref{cf} with $\ti q$ as its central affine curvature. Hence
$$\ti E(x,t,0)\bpm 1\\0\epm= -\frac{1}{k^2} r_{\xi, k}(0) E(x,t,0) r_{-\ti \xi, k}(0)\bpm 1\\0\epm$$
is a solution of \eqref{cf}.  But $\frac{1}{k} r_{\xi,k}(0)=\frac{1}{k}\bpm \xi & \xi^2-k^2\\ 1& \xi\epm$ is a constant in $SL(2,\R)$. So
$$\ti\g= -\frac{1}{k}E(x,t,0) r_{-\ti \xi, k}(0)\bpm 1\\0\epm= -\frac{1}{k}(\g, \g_x) \bpm -\ti \xi\\ 1\epm$$ is a solution of \eqref{cf}.
\end{proof}

\bthm\label{ay} {\rm [BT for \eqref{cf} with $k=0$]}\hfil\par

\ni Let $\g$ be a solution of the flow \eqref{cf}, $q$ its central affine curvature, and $E$ an extended frame for $q$ such that $E(0,0,0)=(\g, \g_x)$ at $(0,0)$.  Write
$$E(x,t,\l)= E_0(x,t)+ E_1(x,t)\l + E_2(x,t)\l^2 + \cdots.$$  Let $\xi\in \R$ be a constant,  and
\beq\label{be} \ti\xi:=\frac{(1,\xi)\g_x}{(1,\xi)\g}. \eeq Here $\g,
\g_x$ are column vectors. Then $E_0(x,t)= E(x,t,0)= (\g, \g_x)$ and
\beq\label{az} \ti \g=e_{12} (-\ti \xi \g+ \g_x) +\bpm
-\ti \xi & \ti \xi^2-k^2\\ 1& -\ti \xi\epm E_1(x,t)\bpm -\ti\xi\\
1\epm \eeq is a solution of \eqref{cf} with central affine curvature
 $\ti q= r_{\xi,0}\ast q= -q+ 2\ti \xi^2$. Moreover,
 $$\ti E(x,t,\l)= \l^{-1} r_{\xi,0}(\l) E(x,t,\l)r_{-\ti \xi(x,t), 0}(\l)$$ is an extended frame for $\ti q$.
\ethm

\begin{proof}
Since $E(\cdot, \cdot, 0)$ and $g=(\g, \g_x)$ are solutions of the linear system \eqref{ac} with the same initial condition, we have  $E(x,t,0)= (\g, \g_x)$.  A simple computation implies that
$$\bpm y_1\\ y_2\epm =E(x,t,0)^{-1}\bpm -\xi\\ 1\epm = (\g, \g_x)^{-1}\bpm -\xi \\ 1\epm= \bpm -(1,\xi)\g_x\\ (1,\xi)\g\epm.$$
Hence $\ti \xi$ defined in  Theorem \ref{ab} is given by \eqref{be}.

By Theorem \ref{ab}, $\ti E(x,t,\l)$ is an extended frame for the new solution $\ti q$. To compute $\ti E(x,t,0)$ we only need to compute the coefficient of $\l$ of $r_{\xi,k}(\l) E(x,t,\l) r_{-\ti\xi(x,t), k}(\l)$. So we get
\beq\label{bf}
\ti E(x,t,0)= e_{12} E_0 r_{-\ti\xi, 0}(0) + r_{\xi,0}(0) E_0 e_{12} + r_{\xi,0}(0) E_1 r_{-\ti \xi, 0} (0).
\eeq
By Proposition \ref{cfkdv}, the first column of $\ti E(x,t,0)$ is a solution of \eqref{cf} with $\ti q$ as its central affine curvature.
 \end{proof}

 As a consequence of Permutability Theorem \ref{bd}, we have
 
 \bthm\label{bo} {\rm [Permutability for \eqref{cf}]\/}\hfil\par

 \ni
 Let $\g$ be a solution of \eqref{cf} with central affine curvature $q$, and $E$ the extended frame for $q$ with $E(x,t,0)= (\g, \g_x)$.  Let $k_i, \xi_i, \ti\xi_i$ be as in Theorem  \ref{bd}, $k_1k_2\not=0$, and $\g_i= \frac{1}{k_i} (\ti\xi_i \g-\g_x)$ for $i=1,2$. Then
 \beq\label{bo1}
 \g_{12}= \frac{1}{k_2}(\ti\xi_{12} \g_1- (\g_1)_x)
 \eeq
  is a solution of \eqref{cf} with central affine curvature $q_{12}= q-2(\ti \xi_1^2 -k_1^2) + 2(\ti \xi_{12}^2-k_2^2)$, where 
  $\ti\xi_{12}= -\ti \xi_1 + (k_1^2-k_2^2)(\ti \xi_1-\ti \xi_2)^{-1}$.
   \ethm

\beg {\rm Explicit solutions} \

We apply B\"acklund transformation (BT) ( Theorem \ref{ax}) to the stationary solution $\g(x,t)= (1, x)^t$ to get explicit solutions of \eqref{cf}, whose central affine curvatures are pure 1-soliton solutions. Then we apply BT again to obtain explicit solutions of \eqref{cf} whose central affine curvature are pure $2$-soliton solutions of the KdV equation.  

First note that
\beq\label{ec}
E(x,t,\l)=\bpm \cosh(zx+z^3t) & z\sinh(zx+z^3t)\\ z^{-1}\sinh(zx+z^3t) & \cosh(zx+z^3t)\epm
\eeq
is an extended frame for the trivial solution $q=0$ of the KdV equation with $\l=z^2$ and
$$\g(x,t)=E(x,t,0)\bpm 1\\ 0\epm= \bpm 1\\ x\epm.$$  
We apply BT to the stationary solution with $k\not=0$ and $\xi=0$ to $\g$. Then we get 
$$\ti \g= \bpm \tanh(kx+k^3t)\\ x\tanh(kx+k^3t)-\frac{1}{k}\epm$$ is a solution of \eqref{cf} with central affine curvature
$$\ti q= -2k^2\sech^2(kx+k^3t).$$

First we use Permutability formula to write down more explicit solutions of \eqref{cf}.  Apply BT to the stationary solution with $\xi=0$ and real $k_1, k_2$ to get two $1$-soliton solutions:
\begin{align*}
&\ti\xi_i=  k_i\tanh (k_ix+ k_i^3t),\\
&\g_i= \bpm \tanh (k_ix+k_i^3t)\\ x\tanh(k_ix+k_i^3t)-\frac{1}{k_i}\epm, 
\end{align*}
with central affine curvature 
$$q_i= -2k_i^2 \sech^2(k_ix+k_it)$$ for $i=1, 2$.  
Then we apply the Permutability formula \eqref{bo1} to get the solution
$$\g_{12}= \frac{1}{k_2} (\ti \xi_{12} \g_1- (\g_1)_x)$$
 of \eqref{cf} with central affine curvature 
$$q_{12}=2k_1^2\sech^2(k_1x+k_1^3t)+2(\ti{\xi}_{12}^2-k_2^2),$$ where
$$
\ti{\xi}_{12}=-k_1\tanh(m_1)+\frac{(k_1^2-k_2^2)(\cosh(m_1+m_1)+\cosh(m_1-m_2))}{(k_1-k_2)\sinh(m_1+m_2)+(k_1+k_2)\sinh(m_1-m_2)}.$$
Note that in general, $q_{12}$ and $\g_{12}$ have singularities. 

\ms
Next we apply BT to $\g_1$ to get new smooth solutions.  Note that 
\beq
E_1(x,t, \l)=\frac{1}{\l^2-k_1^2}\bpm 0 & \l-k_1^2\\ 1 & 0 \epm E(x, t, \l) 
\bpm -\ti{\xi} & \l-k_1^2+\ti{\xi}^2 \\ 1 & -\ti{\xi} \epm\eeq
is an extended frame for $q_1$. 
Apply  BT to $\g_1$ with $k_2\not=0$ and $\xi_2=0$ to get
\begin{align*}
\bpm y_{11} \\ y_{12}\epm & =E(x, t, k_2)^{-1} \bpm 0 \\1\epm \\
&=\bpm \ti{\xi}_1\cosh(k_2x+k_2^3t)-k_2^{-1}(\ti{\xi}_1^2+k_2^2-k_1^2)\sinh(k_2x+k_2^3t) \\
\cosh(k_2x+k_2^3t)-k_2^{-1}\ti{\xi}_1\sinh(k_2x+k_2^3t)\epm,
\end{align*}
where $\ti{\xi}_1=k_1\tanh(k_1x+k_1^3t)$. Let 
$$m_1=k_1x+k_1^3t, \quad m_2=k_2x+k_2^3t.$$ Then we have
$$
\ti{\g}_{12}=-\bpm \frac{\ti{\xi}_1\ti{\xi}_{12}+\ti{\xi}_1^2-k_1^2}{k_1^2}  \\ 
\frac{x(\ti{\xi}_1\ti{\xi}_{12}+\ti{\xi}_1^2-k_1^2)-\ti{\xi}_{12}-\ti{\xi}_1}{k_2^2}\epm.
$$
is a solution of \eqref{cf} with central affine curvature 
$$q_{12}= -q_1 + 2(\ti \xi_{12}^2 -k_2^2),$$ where
$$
\ti{\xi}_{12}=-2\frac{k_1k_2\sinh m_1 \cosh m_1 \cosh m_2-k_2^2\sinh m_2 \cosh^2 m_1+k_1^2\sinh m_2}{\cosh m_1((k_2-k_1)\cosh(m_1+m_2)+(k_1+k_2)\cosh(m_1-m_2))}.
$$
Note that if $k_2>k_1>0$, then $\ti \xi_{12}$ is smooth on $\R^2$. Hence $\g_{12}$ is smooth on $\R^2$ with  smooth $2$-soliton soluiton $q_{12}$ as central affine curvature.  
\eeg

\bigskip


\section{Bi-Hamiltonian structure}\label{ej}

In this section, we first give a brief review of the bi-Hamiltonian structure of the KdV hierarchy (cf. \cite{Mag78}, \cite{GD78}, and \cite{Dic03}) including the sequence of Poisson structures $\{\, ,\}_{2j+1}$ of order $2j+1$. 
Then we use the pull back of $\{\, , \}_{2j+1}$ to get a sequence of compatible Poisson brackets of order $2j-5$ on $\calM(S^1)$ via the map $\Psi: \calM_2(S^1)\to C^\infty(S^1,\R)$ defined in Corollary \ref{bj}. 
We note that the pull back of $\{\, , \}_5$ and $\{\, , \}_3$ on $\calM(S^1)$  give rise to the pre-symplectic forms $\hat w_5$ and $\hat w_3$ given by Pinkall in \cite{UP95} and Fujioka and Kurose in \cite{FK13} respectively.  Moreover, $\hat w_3$ is a symplectic form on $\calM(S^1)/SL(2,\R)$, and $\hat w_5$ is a symplectic form on $\calM(S^1)/S^1$. 

The gradient of a functional $H:C^\infty(S^1, \R)\to \R$ with repect to the $L^2$ inner product  is defined by
$$\rd H_q(v)=\li  \K H(q), v\ri =\oint (\K H(q)) v\rd x.$$
A {\it Poisson operator\/} on $C^\infty(S^1,\R)$ is a collection of linear skew-adjoint operators $L_q$ on $C^\infty(S^1, \R)$ for $q\in C^\infty(S^1,\R)$ such that 
$$\{F_1, F_2\}(q)= -\oint L_q(\K F_1(q))F_2(q)\rd x$$
defines a Poisson structure on $C^\infty(S^1,\R)$.  The Hamiltonian equation for $H:C^\infty(S^1,\R)\to \R$ with respect to $\{\, , \}$ is
$$q_t= L_q (\K H(q)).$$

The following results concerning the bi-Hamiltonian structure of the KdV hierarchy are known:

(i) For $q\in C^\infty(S^1,\R)$, let
 \begin{align}
(L_1)_q(v) &= v_x, \label{dx1}\\
(L_3)_q(v) &= \frac{1}{4}( v_{xxx} - 4qv_x - 2q_xv). \label{dx2}
\end{align}
Then $L_1, L_3$ are Poisson operators on $C^\infty(S^1, \R)$. We let $\{\, ,\}_1$ and $\{\, ,\}_3$ denote the Poisson structures defined by $L_1$ and $L_3$ respectively. 

(ii) Let    $H_{2j+1}:C^\infty(S^1,\R)\to \R$ denote the functional defined by 
\beq\label{eb}
H_{2j+1}(q)= 
\frac{4}{2j+1}\oint C_{-(j+1)}(q)\rd x, \eeq where $Q_i(q)= \bpm
A_i(q) & B_i(q)\\ C_i(q) & - A_i(q)\epm$ is defined in section
\ref{bq}. 
For example,
\begin{align*}
H_1(q)&=-2  \oint q\rd x,\\
H_3(q) &= \frac{1}{2} \oint q^2 \rd x,\\
H_5(q) &= -\frac{1}{8} \oint (q_x)^2 + 2q^3\rd x.
\end{align*}
Then the $(2j+1)$-th flow \eqref{aq} is the Hamiltonian flow for $H_{2j+1}$ ($H_{2j+3}$ resp.) with respect to $\{\, ,\}_3$ ($\{\, ,\}_1$ resp.).
Let 
$$P=L_1^{-1} L_3,$$
where $(L_1)_q^{-1}$ is defined on $\{v\in C^\infty(S^1,\R)\n \oint v\rd x=0\}$ and $$(L_1)_q^{-1}(p(q)_x)= p(q),$$ where $p(q)$ is a polynomial differential of $q$ without constant term. We have
\beq\label{ek}
\K H_{2j+1}(q)= - 2P_q^j(1)= P_q^{j-1}(q), 
\eeq
and the $(2j+1)$-th flow in the KdV hierarchy is 
$$q_{t_{2j+1}}= (L_3)_q(\K H_{2j+1}(q))= (L_1)_q (\K H_{2j+3}(q)).$$ 
  Since the flows in the KdV-hierarchy commute, we have 
$$\{H_{2i+1}, H_{2k+1}\}_1=\{H_{2i+1}, H_{2k+1}\}_3=0, \quad [X_{2i+1}, X_{2k+1}]=0$$ for all $i, k\geq 0$, where 
\beq\label{es}
X_{2i+1}(q)= -2(L_3P^i)_q(1)=(L_3 P^{i-1})_q(q).
\eeq

(iii) The Poisson structures $\{\, , \}_1$ and $\{\, ,\}_3$ are {\it compatible\/}, i.e., $c_1\{\, ,\}_1+ c_3\{\, , \}_3$ is a Poisson structure for all real constants $c_1, c_3$.  This implies that  \beq\label{ei}
 L_{2j+1}= L_3 (L_1^{-1}L_3)^{j-1}, \quad j\geq0,
 \eeq
 is a Poisson operator. We  give a heuristic argument how this sequence $\{L_{2j+1}\}$ arises. 
Since $L_1-\mu L_3$ is a Poisson structure for all $\mu\in \R$, the $2$-form defined by 
$$w(\mu)_q(v_1, v_2)=\li (L_3-\mu L_1)_q^{-1}(v_1), v_2\ri$$ 
is closed.  But 
\begin{align*}
(L_3-\mu L_1)^{-1}&=( L_3(\I- \mu L_3^{-1}L_1))^{-1}= (1-\mu L_3^{-1}L_1)^{-1} L_3^{-1}\\
&= \sum_{i\geq 0} \mu^i (L_3^{-1}L_1)^i L_3^{-1}.
\end{align*}
So we can write 
$w(\mu)= \sum_{i\geq 0} \mu^i w_{2i+3}$,
where 
$$(w_{2i+3})_q(v_1, v_2)= \oint ((L_3^{-1}L_1)^{i-1} L_3^{-1})_q (v_1) )v_2\rd x$$
for $i\geq 0$.  Since $w(\mu)$ is a closed $2$-form for all parameter $\mu$, $w_{2i+3}$ is closed for all $i\geq 0$.
Therefore $L_{2j+3}= L_3P^j$ is an order $(2j+3)$ Poisson operator on the domain of $P^j$ for $j\geq 0$. Note that we can use induction to see that given $v\in C^\infty(S^1\R)$, if $v$  is perpendicular to $X_1(q),\ldots, X_{2j+1}(q)$ then $P^j_q(v)$ is defined (i.e., it is periodic).   

  The $(2j+1)$-th flow in the KdV hierarchy is Hamiltonian flow for $H_{2i+1}$ with respect to $\{\, , \}_{2(j-i)+1}$, i.e.,
 $$q_{t_{2j+1}} = L_{2(j-i)+1}(\K H_{2i+3})= (L_3P^{j-1})_q(q).$$
 Since $L_{2j+1}= L_3 (L_1^{-1}L_3)^{j-1}= (L_3L_1^{-1})^j L_1$, the flows in the KdV hierarchy can be obtained by applying the recursive operator $L_3L_1^{-1}$ as follows:
$$q_{t_{2j+1}}= (L_3L_1^{-1})_q^j (q_x). $$ 

\ms

In the rest of the section, we discuss properties of the pull back of the Poisson structure $\{\, ,\}_{2j+1}$ on the subring 
$$\Psi^*(C^\infty(S^1,\R))=\{F\circ \Psi\n F:C^\infty(S^1,\R)\to \R\}$$ defined by
$$\{F_1\circ \Psi, F_2\circ \Psi\}_{2j+1} = \{F_1, F_2\}_{2j+1}\circ \Psi.$$
In order to write down the corresponding Poisson operator, we need the following Propositions: 

\bprop\label{ed}
Let $\Psi:\calM_2(S^1)\to C^\infty(S^1, \R)$ be the map defined in Corollary \ref{bj}. Then
$$\rd \Psi_\g(\ti \xi) =-2(L_3)_q(\xi),$$
where $q=\Psi(\g)$, $\xi\in C^\infty(S^1, \R)$, $\ti \xi=-\frac{1}{2} \xi_x \g + \xi \g_x$, and $(L_3)_q$ is the Poisson operator defined by \eqref{dx2}.
\eprop

\begin{proof}
Recall that we identify  $C^\infty(S^1, \R)$ as $T(\calM_2(S^1))_\g$ via
$$\xi\mapsto \ti \xi= -\frac{1}{2}\xi_x \g + \xi \g_x.$$
Write $\d\g=\ti\xi$. Take variation of the equation $\g_{xx}= q \g$ to get
\beq\label{dz1}
(\d \g)_{xx}= (\d q)\g + q\d \g = (\d q)\g + q(-\frac{1}{2} \xi_x \g + \xi \g_x)= (\d q -\frac{1}{2} \xi_xq)\g + \xi \g_x.
\eeq
Let $g=(\g, \g_x)$. Then we have
\beq\label{dy}
g_x= g\bpm 0& q\\ 1 &0\epm.
\eeq
Take $x$-derivative of
$$\d \g= -\frac{1}{2}\xi_x \g + \xi \g_x=g \bpm -\frac{\xi_x}{2}\\ \xi\epm $$ and use \eqref{dy} to get
\beq\label{dz2} (\d\g)_{xx}=  (-\frac{1}{2}\xi_{xxx}  + \frac{3}{2}
q\xi_x + q_x \xi)\g + (\xi q ) \g_x. \eeq Compare coefficients of
$\g$ of \eqref{dz1} and \eqref{dz2} to get $\d q=
-\frac{1}{2}\xi_{xxx} + 2 q\xi_x + q_x \xi$, which is equal to
$-2(L_3)_q(\xi)$.
\end{proof}

\bprop\label{eu} Let $C^\infty_I(S^1,\R)$ and $\Psi$ be as in Corollary \ref{eq}, $\g\in \calM_2(S^1)$, and $q=\Psi(\g)$. 
\ben
\item  Let $A\in sl(2,\R)$, and $\hat A(\g)= A\g$ the infinitesimal vector field for the action of $SL(2,\R)$ on $\calM_2(S^1)$.  Then there exists $\xi(\g)\in C^\infty(S^1,\R)$ such that $\hat A(\g)= \widetilde{\xi(\g)}$. In fact, $\xi$ can be computed from $g^{-1}A\g=(-\frac{1}{2}\xi(\g)_x, \xi(\g))^t$, where $g=(\g, \g_x)$.
\item $\Ker((L_3)_q)$ is the space of all $\xi\in C^\infty(S^1,\R)$ such that $\widetilde{\xi(\g)}= A\g$ for some $A\in sl(2,\R)$. 
\item The tangent space $T(C^\infty_I(S^1,\R))_q$ at $q$ is equal to $\Im((L_3)_q)$, and
 $$\Ker (\rd \Psi_\g)=\{\ti\xi\n \xi\in \Ker ((L_3)_q)\}=T(SL(2,\R)\cdot \g)_\g.$$
 \een
\eprop

\begin{proof} It follows from Proposition \ref{ed} and Corollary \ref{bj} that the space of infinitesimal vector fields at $\g$ generated by the $SL(2,\R)$-action on $\calM_2(I)$ is equal to $\Ker((L_3)_q)$. This proves (3).  Since $\hat A(\g)$ is tangent to $\calM_2(S^1)$ at $\g$, there exists $\xi\in C^\infty(S^1,\R)$ such that $\hat A(\g)= \widetilde{\xi(\g)}$.  But 
$$\hat A(\g) = A\g = -\frac{\xi_x}{2} \g + \xi \g_x = (\g, \g_x)\bpm -\frac{\xi_2}{2}\\ \xi\epm = g\bpm -\frac{\xi_2}{2}\\ \xi\epm$$
implies statements (1) and (2). 
\end{proof}

\bprop\label{eh}
Let $H:C^\infty(S^1,\R)\to \R$ be a functional, and $\hat H= H\circ \Psi$.  Then $(\K \hat H)(\g)= \widetilde{y(\g)}$, where
$y(\g)= 2(L_3)_q(\K H(q))$ and $q= \Psi(\g)$
\eprop

\begin{proof}
Use Proposition \ref{ed} to get
$$\d \hat H= \oint (\K H(q))\d q \rd x =\oint -2\K H(q) (L_3)_q(\xi)\rd x=\oint 2 (L_3)_q(\K H(q))\xi \rd x,$$
where $\d \g= \ti \xi$.  
\end{proof}

Next we write down the formula for the Poisson operators induced from $L_{2j+1}$ via the map $\Psi$: 

\bprop\label{ey} The  Poisson operator induced from $L_{2j+1}$ via the map $\Psi$ is 
$$(\hat L_{2j+1})_\g (\ti\xi)= \widetilde{(J_{2j+1})_q(\xi)}, \quad J_{2j+1}= -\frac{1}{4}(L_1^{-1}L_3)^{j-1}L_3^{-1},$$
where $q=\Psi(\g)$.  
\eprop

\begin{proof}
Let $H$ be a functional on $C^\infty(S^1,\R)$, $\hat H= H\circ \Psi$, and $\K\hat H(\g)= \widetilde{y(\g)}$.  
The Hamiltonian vector field of $H$ with respect to $\{\, , \}_{2j+1}$ is $L_{2j+1}(\K H)$.  By Proposition \ref{eh}, $y(\g)= 2(L_3)_q(\K H(q))$. Suppose $\ti\xi$ is the Hamiltonian vector field of $\hat H$ with respect to the Poisson structure induced from $L_{2j+1}$ via $\Psi$.   Then $\rd\Psi(\ti \xi(\g))= L_{2j+1}(\K H)$. By Proposition \ref{ed}, we have
$$\xi(\g)= -\frac{1}{4} L_3^{-1} L_{2j+1} L_3^{-1}(y(\g))= -\frac{1}{4} J_{2j+1}(y(\g)).$$ 
\end{proof}

\brem\ 
\ben  
\item[(a)] The Poisson operator $\hat L_{2j+1}$ has order $(2j-5)$. So the induced Poisson structure loses six derivatives.

\item[(b)] $J_3= -\frac{1}{4}L_3^{-1}$ and $J_5= -\frac{1}{4} L_1^{-1}$.

\item[(c)]  Since $\hat L_3, \hat L_5$ define Poisson brackets, the two forms $\hat w_3, \hat w_5$ defined as below are closed:
\begin{align*}
(\hat w_3)_\g(\ti \xi, \ti \eta)&= \oint ((J_3)_q^{-1}(\xi)) \eta \rd x =-4\oint ((L_3)_q(\xi))\eta\rd x,\\
(\hat w_5)_\g(\ti \xi, \ti \eta)&= \oint ((J_1^{-1})_q(\xi)) \eta\rd x=-4\oint ((L_1)_q(\xi)) \eta\rd x,
 \end{align*}
 where $q=\Psi(\g)$. 
  They are degenerate because $\Ker(L_1)_q=\R$ and $\Ker (L_3)_q$ is of dimension $3$.   The infinitesimal vector field defined by the $S^1$-action on $\calM_2(S^1)$ is $\ti 1(\g)= \g_x$. So $\hat w_5$ is a symplectic form on  $\calM(S^1)/ S^1$. 
It follows from Proposition \ref{eu} that $\hat w_3$ defines a symplectic form on $\calM(S^1)/SL(2,\R)$. 

 \item[(d)] Given $X=\ti \xi=-\frac{\xi_x}{2}\g+ \xi\g_x$ and $Y=\ti \eta= -\frac{\eta_x}{2} \g + \eta\g_x$ in $T\calM_2(S^1)_\g$, we use $\g_{xx}= q\g$ to get
 $$X_x= (-\frac{\xi_{xx}}{2}+ q\xi)\g+ \frac{\xi_x}{2}\g_x, \quad Y_x=  (-\frac{\eta_{xx}}{2}+ q\eta)\g+ \frac{\eta_x}{2}\g_x.$$
A direct computation implies that  
\begin{align}
&(\hat w_3)_\g(X, Y)=-2 \oint \det(X_x,Y_x)+ q\det(X,Y)\rd x,\label{ep1}\\
& (\hat w_5)_\g(X,Y)= -4\oint \det(X,Y)\rd x.\label{ep2}
\end{align}
Note that $\frac{1}{4}\hat w_5$ is the symplectic form \eqref{dp} defined by Pinkall in \cite{UP95}, and $\hat w_3$ is the symplectic form defined by Fujioka and Kurose in \cite{FK13}.

\item[(e)]
The  order $(2j+1)$ central affine curve flow \eqref{cfj} is the Hamiltonian flow for $\hat H_{2(j-i)+1}$ with respect to the Poisson structure $\hat J_{2i+1}$ for all $0\leq i\leq j$. In fact,  \eqref{cfj} is
$$\g_{t_{2j+1}}=-\frac{1}{2}\widetilde{P_q^{j-1}(q)} = \hat L_{2i+1} (\K\hat H_{2(j-i)+3}), \quad 0\leq i\leq j,$$
where $P= L_1^{-1}L_3$.  
\een
\erem


\bigskip

\begin{thebibliography}{99}

\bibitem{AKNS74}Ablowitz, M.J., Kaup, D.J., Newell, A.C., Segur, H., \emph{{T}he inverse scattering transform - Fourier analysis for nonlinear problems}, Stud. Appl. Math \textbf{53} (1974), 249--315.

\bibitem{A81}Adler, M., \emph{On the B\"acklund transformation for the Gel'fand-Dickey equations}, Comm. Math. Phys. \textbf{80(4)} 1981, 517--527.


\bibitem{CIM09}Calini, A., Ivey, T., Mar{\'{\i}} Beffa,G., \emph{Remarks on KdV-type flows on star-shaped curves}, Phys. D \textbf{238(8)} (2009), 788--797.

\bibitem{CIM13}Calini, A., Ivey, T., Mar{\'{\i}} Beffa,G., \emph{Integrable flows for starlike curves in centroaffine space}, SIGMA Symmetry Integrability Geom. Methods Appl. \textbf{9} (2013).

\bibitem{CQ01}Chou,K.S., Qu, C.Z., \emph{The {K}d{V} equation and motion of plane
curves}, J. Phys. Soc. Japan \textbf{70(7)} (2001), 1912--1916.


\bibitem{Dic03}Dickey, L.A., \emph{Soliton equations and Hamiltonian systems}, second edition, Advanced Series in Mathematical Physics \textbf{26} (2003), World Scientific Publishing Co. Inc., River Edge, NJ.

\bibitem{DirSok84}Drinfel'd, V.G., Sokolov, V.V., \emph{Lie algebras and equations of Korteweg-de Vries type},  (Russian) Current problems in mathematics, \textbf{24} (1984),  81--180, Itogi Nauki i Tekhniki, Akad. Nauk SSSR, Vsesoyuz. Inst. Nauchn. i Tekhn. Inform., Moscow.

\bibitem{FK10}Fujioka, A., Kurose, T., \emph{{H}amiltonian formalism for the higher KdV flows on the space of closed complex equicentroaffine curves}, Int. J. Geom. Methods Mod. Phys. \textbf{7(1)} (2010), 165-175.

\bibitem{FK13}Fujioka, A., Kurose, T., \emph{{M}ulti-Hamiltonian structures on space of closed equi-centroaffine plane curves associated to higher KdV flows}, preprint,  Arxiv: math.dg 1310.1688

\bibitem{GD78}Gel'fand, I.M., Dikii, L.A., \emph{{A} family of Hamiltonian structures connected with integrable nonlinear differential equations}, Akad. Nauk SSSR Inst. Prikl. Mat. Preprint \textbf{136} (1978), 41.

\bibitem{GGKM74}Gardner, C.S., Greene, J.M., Kruskal, M.D., Miura, R.M., \emph{{K}orteweg-de Vries equation and generalization. {VI}. Methods for exact solution}, Comm. Pure Appl. Math. \textbf{27} (1974), 97--133.

\bibitem{HS02}Huang, R.P., Singer, D.A., \emph{A new flow on starlike curves in {$\R^3$}}, Proc. Amer. Math. Soc. \textbf{130(9)} (2002), 2725-2735.

\bibitem{L75}Lax, Peter D., \emph{Periodic solutions of the KdV equation}, Comm. Pure Appl. Math. \textbf{28} (1975), 141-188.

\bibitem{Mag78} Magri, F., \emph{A simple model of the integrable Hamiltonian equation}, J. Math. Physics \textbf{19} (1978), 1156-1162.

\bibitem{UP95}Pinkall, U., \emph{Hamiltonian flows on the space of star-shaped
curves},  Results Math. \textbf{27(3-4)} (1995), 328--332.

\bibitem{SS93}Sattinger, D.H., Szmigielski, J.S., \emph{Factorization and the dressing method for the Gel'fand-Dikii hierarchy}, Phys. D \textbf{64(1-3)} (1993), 1--34.


\bibitem{TU00}Terng, C.L., Uhlenbeck, K., \emph{B\"acklund transformations and
loop group actions}, Comm. Pure Appl. Math. \textbf{53} (2000),
1--75.

\bibitem{TU11}Terng, C.L., Uhlenbeck, K., \emph{The $n \times n$ KdV hierarchy}, JFPTA \textbf{10} (2011), 37--61.

\bibitem{TW2} Terng, C.L., Wu, Z., \emph{Central affine curve flows on $\R^n\setminus \{0\}$}, preprint

\bibitem{ZF71}Zaharov, V.E., Faddeev, L.D., \emph{The Korteweg-de Vries equation is a fully integrable Hamiltonian system}, Funkcional. Anal. Priloz \textbf{5(4)} (1971), 18--27.

\end{thebibliography}
\end{document}